\documentclass[12pt]{article}
\usepackage{mathrsfs}
\usepackage{amsfonts}
\usepackage{amssymb}
\usepackage{enumerate}
\usepackage{amsmath}
\usepackage{graphicx}
\usepackage{hyperref}
\usepackage{color}

= 0pt \leftmargin     = 1.25in \topmargin =0pt \headheight     = 0pt \headsep
= 0pt \topskip =0pt
\def\sqr#1#2{{\vcenter{\vbox{\hrule height.#2pt
              \hbox{\vrule width.#2pt height#1pt \kern#1pt \vrule width.#2pt}
              \hrule height.#2pt}}}}
\def\signed #1{{\unskip\nobreak\hfil\penalty50
              \hskip2em\hbox{}\nobreak\hfil#1
              \parfillskip=0pt \finalhyphendemerits=0 \par}}
\def\endpf{\signed {$\sqr69$}}

\def\3n{\negthinspace \negthinspace \negthinspace }
\def\2n{\negthinspace \negthinspace }
\def\1n{\negthinspace }

\def\={\buildrel \triangle \over =}

\def\no{\noindent}

\def\bs{\bigskip}

\def\|{\Big |}
\def\({\Big (}
\def\){\Big )}
\def\[{\Big[}
\def\]{\Big]}

\newtheorem{remark}{Remark}[section]

\newtheorem{theorem}{Theorem}[section]

\newtheorem{definition}{Definition}[section]

\newtheorem{proposition}{Proposition}[section]

\makeatletter
   
   \@addtoreset{equation}{section}
\makeatother
\begin{document}
	\title{\bf  Algebraic structure of the Gaussian-PDMF space and applications on fuzzy equations}
	\author{ 		Chuang Zheng\thanks{ School of Mathematical Sciences,
			Beijing Normal University,
			100875 Beijing,
			China.   Email: chuang.zheng@bnu.edu.cn.
	}}
	%\date{}
	\maketitle 
	
	\begin{abstract}
		In this paper, we extend the research presented in \cite{2023WangZheng} by establishing the algebraic structure of the Gaussian Probability Density Membership Function (Gaussian-PDMF) space.	We consider fixed objective and subjective entities, denoted as $(h,p)$, and provide the explicit form of the membership function. Consequently, every fuzzy number with the membership function in $X_{h,p}(\mathbb{R})$, denoted as $\tilde{x}$, can be uniquely identified by a vector $\langle x; d^-, d^+, \mu^-,\mu^+\rangle$. Here, $x\in \mathbb{R}$ represents the ``leading factor" of the fuzzy number $\tilde{x}$ with a membership degree equal to $1$. The parameters $d^-$ (left side) and $d^+$ (right side) denote the lengths of the compact support, while $\mu^-$ (left side) and $\mu^+$ (right side) represent the shapes.  We introduce five operators: addition, subtraction, multiplication, scalar multiplication, and division. We demonstrate that, based on our definitions, the Gaussian-PDMF space exhibits a well-defined algebraic structure. For instance, $X_{h,p}(\mathbb{R})$ is a vector space over $\mathbb{R}$, featuring a subspace that forms a division ring, allowing for the representation of fuzzy polynomials, among other properties. We provide several examples to illustrate our theoretical results.
	\end{abstract}
	
	\bs

	\no{\bf Key Words}. Fuzzy numbers; Gaussian probability density membership function (Gaussian-PDMF) space; Fuzzy arithmetic; Abelian group; Ring; Fuzzy Vieta's formula; Linear space; Fuzzy equations.
	
	\section{Introduction}
	
	Fuzzy numbers and fuzzy set theory are topics originated from Zadeh (\cite{Zadeh1965}) by dealing with the imprecise quantities and uncertainty. Since then, they have found successful applications in a wide range of areas, including pure and applied mathematics, computer science, and related fields such as fuzzy logic, fuzzy information, soft computing, and fuzzy control.
	
	In general,  a fuzzy number can be uniquely determined by its membership function. The construction of  an appropriate membership function is the cornerstone upon which fuzzy set theory has evolved. For existing references,  see \cite{dombi1990, krusinska1986note, mashchenko2021sums, 1996Course, wu2013decomposition, wu2020arithmetic}  and the references therein.  Particularly, in the realm of artificial intelligence, linguistic variables are associated with a membership function based on subjective choices (see, for instance, \cite{Costa2021closure,Li2005fuzzy}). 
	
	Once membership functions of fuzzy numbers are identified, a fundamental concern is performing arithmetic operations on them. In the fuzzy world, arithmetic operations on real numbers at the classical crisp set level translate into algebraic operations on membership functions. Numerous papers have explored the fundamental definitions of these operations and the corresponding algebraic structures. Triangular norms, for example, address binary operations on the interval $[0,1]$ (\cite{KLEMENT20045}). Interactive fuzzy numbers and their arithmetic operators are proposed through joint possibility distributions (\cite{Zadeh1975199, Dubois1981TAC, Esmi2019246}), with a notable method offering minimum norm via the sup-J extension principle (\cite{ESMI2021}). Intuitionistic fuzzy sets, along with probabilistic addition, have successfully addressed fuzzy aggregation problems in expert systems (\cite{ATANASSOV198687,MR1718470, SZMIDT2000505, xu2007}).
	Moreover, in the context of fuzzy differential equations, a common definition of (generalized) Hukuhara differentiability necessitates designing the difference between two fuzzy numbers using $\alpha$-cuts of their corresponding membership functions (\cite{dubois1982towards, hukuhara1967, puri1983differentials}).  Applications encompass various classes of membership functions such as triangle MF, trapezoidal MF, Sigmoid MF, among others. MATLAB's standard toolbox offers more than 11 different types of membership functions for constructing fuzzy models.

	In \cite{2023WangZheng}, a new class of membership functions is introduced through the composition of two types of nonlinear functions: $h:[0,1] \rightarrow \mathbb{R}$ and $p:\mathbb{R}\rightarrow [0,1]$. Here, $h$ represents subjective perception, while $p$ provides objective information via the probability density function. This class of functions defines the Probability Density Membership Function (PDMF) space $X_{h,p}(\mathbb{R})$. Specifically, $h$ is chosen as the tangent function, and $p$ is given by the Gaussian kernel $p(\cdot, \mu)$ with $\sigma=1$ and $\mu$ undetermined, where $\mu$ measures fuzziness, and arithmetic operations are designed via $\mu$ to achieve a linear structure. Such class of MFs is called as the Gaussian-PDMFs. Notably, the $\alpha$-cuts representation is unnecessary during computation implementation.

	In this paper, we introduce a new notation of the fuzzy number $\tilde{x}$ with the membership function in Gaussian-PDMF space $X_{h,p}(\mathbb{R})$. As we show later, $\tilde{x}$ is uniquely identified by a vector $[x; d^-, d^+, \mu^-,\mu^+]$, where $x\in \mathbb{R}$ represents the ``leading factor” of the fuzzy number $\tilde{x}$ with membership degree equaling to $1$, $d^-$ (resp. $d^-$) is the length of the compact support on the left side (resp. right side) and $\mu^-$ (resp. $\mu^+$) represents the shape of the left (resp.  right) side.  We design five algebraic operations: addition, subtraction, multiplication, division and scalar multiplication in terms of the five parameters of $\tilde{x}$. 
	
	The first result is that addition $\oplus$ forms an abelian group on $X_{h,p}(\mathbb{R})$. Due to the existence of the inverse element, the minus $\ominus$ can be seen, or directly defined by, the addition of the inverse element. Furthermore,  together with the multiplication $\otimes$,  $X_{h,p}(\mathbb{R})$ turns to a commutative ring with identity. Following the standard procedure of the algebra, we state several basic properties of the space $X_{h,p}(\mathbb{R})$ with two operations $\oplus$ and $\otimes$. Collecting those elements with inverse under multiplication, we form a subset of $X_{h,p}(\mathbb{R})$, denoted by $U(X)$. As a result, the cancellation law holds with respect to $\oplus$ and $\otimes$. We give several applications for this algebraic structure, such as the solvability of fuzzy equations, Binomial Theorem, fuzzy Vieta's formula, etc. The last scalar multiplication on $X_{h,p}(\mathbb{R})$  is the basic operation for the linear algebra. By constructing a basis of $X_{h,p}(\mathbb{R})$, we give the coordinates of $\tilde{x}$ and shows that the coordinate mapping is a isomorphism from  $X_{h,p}(\mathbb{R})$ onto $\mathbb{R}^5$.

	The structure of the paper is organized as follows. In Section \ref{sec2}, we present the basic concepts of fuzzy numbers and outline the requirements of the membership function. In Section \ref{sec3}, we provide exact definitions of arithmetic operators on $X_{h,p}(\mathbb{R})$, along with the algebraic structure such as abelian group, commutative ring with identity, etc. Additionally, we introduce a binomial theorem and demonstrate its application to second-order fuzzy equations. In Section \ref{sec4}, we prove that $X_{h,p}(\mathbb{R})$ is a vector space and establish the isomorphism between $X_{h,p}(\mathbb{R})$ and $\mathbb{R}^5$. In Section \ref{sec5}, we present numerical examples and corresponding graphs to illustrate the operations under consideration. Finally, in Section \ref{sec6}, we offer a concluding remark to provide a comprehensive summary of the paper.

\section{Preliminary}\label{sec2}

The formulation of membership functions is the crucial step in the design of fuzzy system. There are several methods to develop them. We summarize some of them as follows:
\begin{enumerate}[1)]
	\item L-R linear functions, which is the simplest possible model (\cite{dubois1980fuzzy}); 
	\item Rational functions of polynomials (\cite{giachetti1997parametric, guerra2005approximate});
	\item B-Spline MF (\cite{wang1995fuzzy}); 
	\item Piecewise linear functions (\cite{wen2019novel} and refs [8-20] in it);
	\item Gaussian PDMF and the corresponding Gaussian-PDMF space $X_{h,p}(\mathbb{R})$ (\cite{2023WangZheng}).
\end{enumerate}
In all of these definitions of fuzzy numbers, the membership function needs to satisfy the following assumptions:
\begin{definition}\label{1st_Def}% (\cite{SHEN2020})
	A fuzzy number $A$ is a fuzzy subset of the real line $\mathbb{R}$ with membership function $f_{A}$ which possesses  the following properties: 
	\begin{enumerate}[a)]
		\item $f_{A}$ is fuzzy convex,
		\item $f_{A}$ is normal i.e., $ \exists x_{0} \in \mathbb{R} $ such that $f_{A}(x_{0})=1$,
		\item $f_{A}$ is upper semi-continuous,
		\item The closure of the set $\{x\in\mathbb{R}|f_{A}(x) > 0\}$ is compact.
	\end{enumerate} 
\end{definition}

Definition \ref{1st_Def} is  straightforward and has been used extensively in practical applications  (\cite{SHEN2020}). However, the above conditions are too vague  and a particular class of functions, named as monotonic fuzzy numbers, is introduced  with the  following more precise assumptions:

\begin{definition}\label{2nd_Def}
	A monotonic fuzzy number $\tilde{b}$, denoted by $\tilde{b}=(a,b,c)$,  is defined as a membership function $f(x)$  which possesses  the following properties (\cite{dubois1980fuzzy}): 
	\begin{enumerate}[a)]
		\item $f(x)$ is increasing on the interval $[a,b]$ and decreasing on $[b,c]$,
		\item $f(x)=1$ for $x=b$, $f(x)=0$ for $x \leq a$ or $x \geq c $,
		\item $f (x)$ is upper semi-continuous,
	\end{enumerate} 
where $a,b,c,$ are real numbers satisfying $-\infty < a< b < c < +\infty$. 
\end{definition}  
Clearly, a class of triangular fuzzy numbers is a subset of the class of monotonic fuzzy numbers. It is due to the fact that, in the definition of the triangular fuzzy numbers, the function in the condition $a)$ of Definition \ref{2nd_Def} is restricted by linear ones  (\cite{liang2013}).

Recently,  a new class of fuzzy numbers called PDMF is given in \cite{2023WangZheng}, in which extra information between $a,b$ and $b,c$ are clarified. More precisely, a monotonic number $\tilde{b}$, denoted by $\tilde{b}=\langle (a,b,c);P(x^-, y^-), Q(x^+,y^+)\rangle $,  is defined as a membership function $f(x)$  which possesses  the following structure:

\begin{equation}\label{fxabc}
	\begin{array}{lll}	
	f(x)	
	&=&\left \{
	\begin{aligned}
	&0, &x \in (-\infty,a]\\
	&\int_{-\infty}^{h^{-}(x)}p^{-}(y)dy,  &x\in (a,b)\\
	& 1 , &x=b\\
	&\int_{-\infty}^{h^{+}(x)}p^{+}(y)dy,  &x\in (b,c)\\
	&0,    &x\in [c,+\infty).\\
	\end{aligned}
	\right.
	\end{array}
	\end{equation} 

In \eqref{fxabc}, $h^{-}, h^{+}$ are the auxiliary functions satisfying 
\begin{enumerate}[i)]
	\item $\lim_{x \rightarrow a^+}h^{-}(x)=\lim_{x \rightarrow c^-}h^{+}(x)=-\infty, \quad \lim_{x \rightarrow b^-}h^{-}(x)=\lim_{x \rightarrow b^+}h^{+}(x)=+\infty$,
	\item $h^{-}$ is continuous and increasing on $ (a,b)$,  $h^{+}$ is continuous and decreasing on $ (b,c)$,
\end{enumerate} 
and $p^{-},p^{+}$ are probability density functions (PDFs) satisfying
$$
\displaystyle\int_{-\infty}^{+\infty}p^-(y)dy=1,\;\int_{-\infty}^{+\infty}p^+(y)dy=1, \;
 p^-(t) \geq 0, \; p^+(t) \geq 0,\quad\forall t \in (-\infty,+\infty).
 $$·
 Moreover, it is required that the shape of $f(x)$ passes through $P, Q$.

The following Theorem stated in \cite{2023WangZheng} provides the well-posedness of the definition \eqref{fxabc} for $f(x)$:

\begin{theorem}\label{PDMF} 
	Let 
	$$
	h^{-}(x)=h(\frac{x-a}{b-a}), 
h^{+}(x)=h(\frac{c-x}{c-b})\quad \hbox{with}\quad \lim_{x \rightarrow 0^+}h (x)=-\infty,   \lim_{x \rightarrow 1^-}h (x)=+\infty,
	$$
	and $h(x)$ is continuous and increasing on $ (0,1)$. Let $p^-$ and $p^+$ be originated from the same class $p$ of PDFs.
	Then	there exists  at least one pair $(h, p)$ such that the graph of $f_{h,p}$ passes through $P$ and $Q$, i.e., $f_{h,p}(x^-)=y^-$ and $f_{h,p}(x^+)=y^+$.
	Moreover, the PDMF in the space  
	\begin{equation*}%\label{PDMFS}
		X_{h,p}(\mathbb{R}):=\{f_{h,p}(x):\mathbb{R} \rightarrow [0,1]  \;\; \hbox{is as the form of \eqref{fxabc}} :a < b < c\}
	\end{equation*}
	fulfills all requirements in Definition \ref{1st_Def} and Definition \ref{2nd_Def}. 
	\end{theorem}

\begin{figure}[htbp]\label{Fig.1}
	\centering
	\includegraphics[height=6cm,width=8cm]{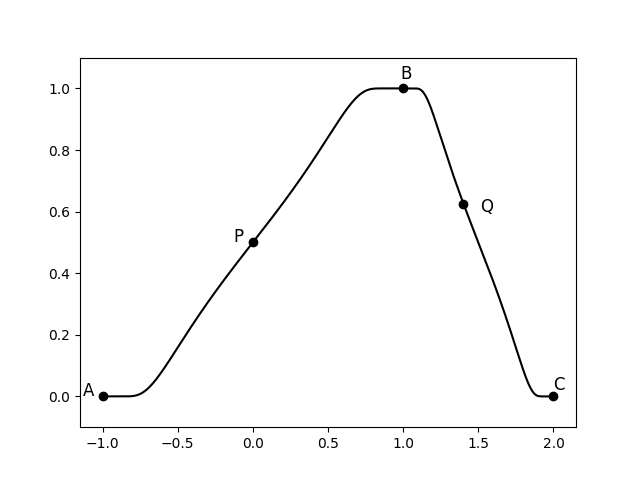}\caption{ Fuzzy number $\tilde{1}=\langle(-1,1,2);P,Q\rangle$ or $\langle 2; 2,1,\mu^{-},\mu^{+}\rangle$} where $\mu^{-}$ is uniquely identified by the point $P$ via $f_-(x)$ in \eqref{mu-mu+}.
\end{figure}

In this paper, to specify the class of functions under consideration, we fix $h$ as the tangent function and $p$ is given by the Gaussian kernel $p(\cdot, \mu)$ with $\sigma=1$ and $\mu$ to be undetermined, i.e. 
	\begin{equation}\label{hp}
	h(x) =\tan(\pi x - \frac{\pi}{2}),  \; x \in (0,1), 
	\qquad p(t)=p (t; \mu)=\frac{1}{\sqrt{ 2\pi}}e^{-\frac12(t-\mu)^{2}}, \; t\in\mathbb{R}.
	\end{equation} 
	The two control points $P(x^{-}, y^{-}), Q(x^{+}, y^{+})$ are placed on each side of the central point $B(x_0, 1)$, and the shape of the function passes through $A(x_0-d^-,0)$ and $C(x_0+d^+,0)$ with $d^-=x_0-a$ and $d^+=c-x_0$, respectively. See Figure \ref{Fig.1} for the sample of the shape of the membership function. Note that here $d^-$ and $d^+$ are the deviation measures from the center $x_0$ to the left side and right side, respectively. 
	Clearly, taking into account \eqref{hp}, the function in the PDMF space has the explicit form 
\begin{equation}\label{abcmumu}
f_{h,p}(x)=
\left \{
\begin{aligned}	
&0, &x \in (-\infty,x_0-d^-]\\
&f_{-}(x;x_0,d^-,\mu^{-}),  &x\in (x_0-d^-,x_0)\\
& 1 , &x=x_0\\
&f_{+}(x;x_0,d^+,\mu^{+}),  &x\in (x_0,x_0+d^+)\\
&0,    &x\in [x_0+d^+,+\infty)\\
\end{aligned}
\right.
\end{equation}
where $f_-$ and $f_+$ are given by
\begin{equation}\label{mu-mu+}
\begin{aligned}
	f_{-}(x;x_0,d^-,\mu^{-})
	&=\int_{-\infty}^{\tan(\frac{\pi}{d^-}(x-x_0+d^-)-\frac{\pi}{2})}\frac{1}{\sqrt{ 2\pi}}e^{-\frac12(t-\mu^{-})^{2}}dt, &\qquad x\in(x_0-d^-,x_0),\\ 
	f_{+}(x;x_0,d^+,\mu^{+})
		&=\int_{-\infty}^{\tan(\frac{\pi}{d^+}(x_0+d^+-x)-\frac{\pi}{2})}\frac{1}{\sqrt{ 2\pi}}e^{-\frac12(t-\mu^{+})^{2}}dt, &\qquad x\in(x_0,x_0+d^+).\\ 
	\end{aligned}
\end{equation}
The function  as above is the {\bf {Gaussian-PDMF}} and the corresponding function space is the {\bf {Gaussian-PDMF Space}}.  More precisely,  
\begin{equation}\label{PDMFS}
	X_{h,p}(\mathbb{R}):=\{f_{\tan,G}(x):\mathbb{R} \rightarrow [0,1]  \;\; \hbox{is as the form of 
	\eqref{abcmumu} with \eqref{mu-mu+}} \}.
    %\eqref{abcmumu}} \}. %: d^->0,d^+>0\}.
\end{equation} 

As described in Theorem $4.1$ of \cite{2023WangZheng}, there is a unique vector $\langle (x_0-d^-,x_0,x_0+d^+);\mu^-,\mu^+\rangle$ for each membership function $f_{\tan,G}(x)$ in the Gaussian-PDMF space passing through the two control points $P(x^{-}, y^{-}), Q(x^{+}, y^{+})$. 

Consequently, it is reasonable to give two equivalent notations of the Gaussian-PDMF as  
\begin{equation}\label{abcde}
 \langle(x_0-d^-,x_0,x_0+d^+);P,Q\rangle
\Longleftrightarrow
\langle x_0; d^-, d^+,\mu^{-},\mu^{+}\rangle. 
\end{equation}
 Here, $(\mu^{-},\mu^{+})$ can be uniquely identified by means of the inverse functions of Formula \eqref{mu-mu+} via the information on two control points $P(x^{-}, y^{-}), Q(x^{+}, y^{+})$.  Note that in the classical Zadeh's summation notation  (see, for instance, Chapter $2$ of \cite{1996Course}), the above information gives $ \tilde{x}_0=0/(x_0-d^{-})+y^-/x^-+1/x_0+y^+/x^++0/(x_0+d^{+})$.

 In the sequel, we will use the new notation $\langle x_0; d^-, d^+,\mu^{-},\mu^{+}\rangle$ for the fuzzy number $\tilde{x}_0$ with the membership function \eqref{abcmumu}-\eqref{mu-mu+}. The advantage of the revision is that the assumption $a<b<c$ is replaced by the positivity of $d^-$ and $d^+$. We emphasize that in \eqref{abcde}, $x_0$ represents the crisp part of the fuzzy number $\tilde{x}_0$ and the other four parameters represent its fuzzy part, respectively.  More precisely, the fuzzy number $\tilde{x}_0$ can be identified by three parts: Value $x_0$, Ambiguity $d^-, d^+$ and Fuzziness $\mu^{-},\mu^{+}$ (see \cite{delgado1998,delgado1998canonical,Rahman2015} and references therein). They capture the relevant information contained in the fuzzy number and could be useful for developing ranking procedure of fuzzy numbers. For instance, Value $x_0$ has the key property of permitting us to associate a real value with the fuzzy number $\tilde{x}_0$. This approach is beyond the scope of this paper and  will  be considered elsewhere.

\section{Arithmetic operators on $X_{h,p}(\mathbb{R})$}\label{sec3}

We now define the addition, multiplication, subtraction and  scalar multiplication on the Gaussian-PDMF space $X_{h,p}(\mathbb{R})$ in \eqref{PDMFS}.
% \begin{equation}\label{GPDMFS}
% 	X_{h,p}(\mathbb{R})
% 	=
% 	\{\tilde x\|\tilde x=\langle x; d^-, d^+,\mu^-,\mu^+ \rangle\;\hbox{has the form of \eqref{abcmumu}} \}.
% 	\end{equation}
For convenience we shall generally use the notation $\tilde x_i $ instead of $\tilde x_i=\langle x_i; d^-_i, d^+_i,\mu^-_i,\mu^+_i\rangle$,  for any $i=0,1,2,\cdots$.
	\begin{definition} \label{def} Let $\tilde{x}_{1}, \tilde{x}_{2}$ be two Gaussian-PDMFs in $X_{h,p}(\mathbb{R})$, we define
		\begin{enumerate}[(1)]
		\item  $\tilde{x}_{1}\oplus\tilde{x}_{2}=\langle x_1+x_2; d^-_1d^-_2, d^+_1d^+_2, \mu^{-}_{1}+\mu^{-}_{2},\mu^{+}_{1}+\mu^{+}_{2} \rangle$,
		\item 
		$\lambda \tilde{x}=
		\langle \lambda x; (d^-)^\lambda,(d^+)^\lambda, \lambda \mu^{-},\lambda\mu^{+} \rangle,
		$
		\item
		$\tilde{x}_{1}\ominus\tilde{x}_{2}=\langle x_{1}-x_{2}; d^-_1 (d^-_2)^{-1}, d^+_1(d^+_2)^{-1},\mu^{-}_{1}-\mu^{-}_{2},\mu^{+}_{1}-\mu^{+}_{2} \rangle$,
		\item  
		$\tilde{x}_{1}\otimes\tilde{x}_{2}=\langle x_{1}x_{2}; (d^-_1)^{\ln d^-_2}, (d^+_1)^{\ln d^+_2},\mu^{-}_{1}\mu^{-}_{2},\mu^{+}_{1}\mu^{+}_{2} \rangle$.
	\end{enumerate}
	\end{definition}

	Some remarks are in order.

	\begin{remark}\label{crisp number}
		Note that $a<b<c$ implies  $d^-,d^+\neq 0$. Hence the crisp number $x\in\mathbb{R}$ is excluded in the Gaussian-PDMF space since it coincides to  the class of vectors $\langle x;0,0,\cdot,\cdot\rangle$ where $d^-=d^+=0$.  In this case, each  crisp number is formally equivalent to a class 
		$$
	\{ \tilde{y}=\langle x; 0, 0,\mu^{-}_y,\mu^{+}_y,\rangle, \;\forall \mu^{-}_y,\mu^{+}_y\in\mathbb{R}\}.
		$$
	Consequently, $x\in\mathbb{R}$ can be seen as the ``limitation" of some quotient class of $X_{h,p}(\mathbb{R})$ satisfying 
		$$
		 x \sim \lim_{d^-,d^+\rightarrow 0}\{\tilde{x}\in X_{h,p}(\mathbb{R}) \;\hbox{for any }\; \mu^{-}_y,\mu^{+}_y\in\mathbb{R}\}.
		$$ 
		and all arithmetic operations on the Gaussian-PDMF space will coincide with the corresponding operations on $\mathbb{R}$.
\end{remark}
\begin{remark}\label{ranking}
		Note that in the existing literature, ranking fuzzy numbers are  one of the most important research topics in fuzzy set theory. 	Mathematically speaking,  ranking fuzzy numbers can be seen as ranking the functions as the form of \eqref{abcmumu} in Gaussian-PDMF space. 
		One simple ranking method to evaluate Gaussian-PDMF numbers is that: for any $\tilde{x}, \tilde{y}\in X_{h,p}(\mathbb{R})$, $\tilde{x}\leq \tilde{y}$ if and only if $x\leq y$.  As shown later, all properties of crisp numbers continues to keep under Definition \ref{def}.
\end{remark}
\begin{remark}
		One of the main advantage of Definition \ref{def} is that, as we show later, the nice algebraic structure of Gaussian-PDMF makes the computation to be universal. For instance, it is no need to add the assumptions that the fuzzy numbers involved have to has same signs, nor the existence of the fuzzy number of H/gH difference. Furthermore, it is not necessary to use the $\alpha$-cuts of the fuzzy numbers to perform calculations.
\end{remark}

The first Theorem on the algebraic structure of Gaussian-PDMF space $X_{h,p}(\mathbb{R})$ under Definition \ref{def} is:
\begin{theorem}\label{abeliangroup}
	$X_{h,p}(\mathbb{R})$ together with the binary operation $\oplus$ is an {\bf abelian group}. More precisely, 
	\begin{enumerate}[(i)]
		\item  The binary operation $\oplus$ is {\bf associative}:
		\begin{equation}
			(\tilde{x}_{1}\oplus\tilde{x}_{2})\oplus\tilde{x}_{3}
			=
			\tilde{x}_{1}\oplus(\tilde{x}_{2}\oplus\tilde{x}_{3}), \quad\hbox{for all }\; \tilde{x}_{i}\in X_{h,p}(\mathbb{R}),\;\; i=1,2,3.
		\end{equation}
		\item There exists a unique {\bf  two-sided identity element} $\tilde 0=\langle 0;1,1,0,0\rangle $ such that 
		\begin{equation}
			\tilde 0\oplus \tilde x= \tilde x\oplus\tilde 0, \quad\hbox{for all }\; \tilde{x}\in X_{h,p}(\mathbb{R}).
		\end{equation}
		\item For every $\tilde{x}\in X_{h,p}(\mathbb{R})$  there exists a {\bf two-sided inverse element} $- \tilde{x}\in X_{h,p}(\mathbb{R})$ such that 
		\begin{equation}
			- \tilde{x} \oplus\tilde{x}=   \tilde{x} \oplus - \tilde{x}= \tilde 0. 
		\end{equation}
		Moreover, $- \tilde x=\langle -x;(d^-)^{-1},(d^+)^{-1}, - \mu^{-},-\mu^{+}\rangle  $ coincides to the scalar multiplication $\lambda \tilde{x}$ with $\lambda=-1$.
		\item The binary operation $\oplus$ is {\bf  commutative}:
		\begin{equation}
			\tilde{x}_{1}\oplus\tilde{x}_{2}=\tilde{x}_{2}\oplus\tilde{x}_{1} , \;\hbox{for all }\; \tilde{x}_{i}\in X_{h,p}(\mathbb{R}),\;\; i=1,2.
		\end{equation}
	\end{enumerate}
\end{theorem}
\begin{remark}
	Compare to the Hukuhara difference and generalized-Hukuhara difference (see, for instance, \cite{hukuhara1967,stefanini2010} and references therein), the difference $\ominus$ in Definition \ref{def} shows a clear and simple algebraic structure of the fuzzy numbers. Hence, the subtraction of $\tilde{x}$ and the addition of its negative $-\tilde{x}$ are compatible, which is not true for most of the existing definitions.
\end{remark}

{\bf Proof of Theorem \ref{abeliangroup}:}  For $(i)$, set $\tilde{x}_{i}=\langle x_i; d_i^-,d_i^+,\mu_i^-,\mu_i^+\rangle, i=1,2,3$.  Recall Formula $(1)$ of Definition \ref{def}, it holds
	\begin{equation*}
		\begin{array}{lll}
			(\tilde{x}_{1}\oplus\tilde{x}_{2})\oplus\tilde{x}_{3} 
			& = & \langle x_1+x_2; d_1^-d_2^-,d_1^+d_2^+r,\mu_1^-+\mu_2^-,\mu_1^++\mu_2^+\rangle \oplus \langle x_3; d_3^-,d_3^+,\mu_3^-,\mu_3^+\rangle\\
			& = & \left\langle  x_1+x_2+x_3; d_1^-d_2^-d_3^-, d_1^+d_2^+d_3^+, \mu_1^-+\mu_2^-+\mu_3^-, \mu_1^++\mu_2^++\mu_3^+\right\rangle \\
            & = & \left\langle  x_1+(x_2+x_3); d_1^-(d_2^-d_3^-), d_1^+(d_2^+d_3^+), \mu_1^-+(\mu_2^-+\mu_3^-), \mu_1^++(\mu_2^++\mu_3^+)\right\rangle \\
			& = & \tilde{x}_{1}\oplus(\tilde{x}_{2}\oplus\tilde{x}_{3}).
		\end{array}
    \end{equation*}
	 
    For $(ii)$, it is trivial that $\tilde 0=\langle 0;1,1,0,0\rangle $ is a two-sided identity. If $\tilde 0'$ is also a two-sided identity, then $\tilde 0 = \tilde 0 \oplus \tilde 0'$=$\tilde 0'$, which means that $\tilde 0$ is unique.
    
    Assertion $(iii)$ is straightforward since for any $d^-, d^+\neq 0$, it holds
    $$
    - \tilde{x} \oplus\tilde{x}
    =
    \langle -x;(d^-)^{-1},(d^+)^{-1}, - \mu^{-},-\mu^{+}\rangle    \oplus \langle x; d^-,d^+,  \mu^{-},\mu^{+}\rangle = \langle 0;1,1,0,0\rangle.
    $$
    Clearly, due to Formula $(2)$ of Definition \ref{def}, $(-1)\tilde{x} = -\tilde{x}$.

    Finally, the commutativity of $\oplus$ is a direct consequence of Formula $(1)$ of Definition \ref{def} due to the fact that  $d^-_1d^-_2=d^-_2d^-_1$ and $d^+_1d^+_2=d^+_2d^+_1$.
    \endpf
		As a direct consequence of Theorem \ref{abeliangroup}, the linear fuzzy equation can be solved with the following formula:
		\begin{proposition}\label{linearequation1}
			Let $\tilde{b},\tilde{c}\in X_{h,p}(\mathbb{R})$ with $\tilde{b}=\langle b; d_b^-,d_b^+,\mu_b^-,\mu_b^+\rangle$ and $\tilde{c}=\langle c; d_c^-,d_c^+,\mu_c^-,\mu_c^+\rangle$. Let $a\in\mathbb{R}\setminus \{0\}$. Set $\tilde{x}=\langle x; d^-,d^+,\mu^-,\mu^+\rangle$ be an undetermined fuzzy number such that 
		\begin{equation}\label{1order}
			a\tilde{x}\oplus\tilde{b}=\tilde{c}.
		\end{equation}
		Then there exists a unique fuzzy number $\tilde{x}$ belonging to $X_{h,p}(\mathbb{R})$ satisfying \eqref{1order} with the formula 
		\begin{equation}\label{Soln1}
			\tilde{x}=\frac{1}{a}(\tilde{c}\ominus\tilde{b})=
			\left\langle \frac{c-b}{a}; \left(\frac{d_c^-}{d_b^-}\right)^{\frac{1}{a}},\left(\frac{d_c^+}{d_b^+}\right)^{\frac{1}{a}},
			\frac{\mu_c^--\mu_b^-}{a},\frac{\mu_c^--\mu_b^-}{a}\right\rangle.
		\end{equation}
		\end{proposition}

		In the next theorem, we show that $X_{h,p}(\mathbb{R})$ is a commutative ring together with the operations addition and multiplication.

		\begin{theorem}\label{ring}
			$X_{h,p}(\mathbb{R})$ together with  $\oplus$ and $\otimes$ is a {\bf commutative ring with identity}  $\tilde{1}_X=\langle 1; e,e,1,1\rangle$. More precisely, 
			\begin{enumerate}[(i)]
				\item  $(X_{h,p}(\mathbb{R}),\oplus)$ is an abelian group;
				\item  $(\tilde{x}_1\otimes\tilde{x}_2)\otimes\tilde{x}_3=
				 \tilde{x}_1\otimes(\tilde{x}_2\otimes\tilde{x}_3)$ for all $\tilde{x}_{i}\in X_{h,p}(\mathbb{R}),\; i=1,2,3$ ({\bf associative multiplication});
				 \item $\tilde{x}_1\otimes(\tilde{x}_2\oplus\tilde{x}_3)=
				 \tilde{x}_1\otimes\tilde{x}_2\oplus\tilde{x}_1\otimes\tilde{x}_3$ and $(\tilde{x}_1\oplus \tilde{x}_2)\otimes\tilde{x}_3=
				 \tilde{x}_1\otimes\tilde{x}_3\oplus\tilde{x}_2\otimes\tilde{x}_3$ for all $\tilde{x}_{i}\in X_{h,p}(\mathbb{R}),\; i=1,2,3$ ({\bf left and right distributive laws});
				 \item $\tilde{x}_1\otimes\tilde{x}_2=\tilde{x}_2\otimes\tilde{x}_1$ for all $\tilde{x}_{i}\in X_{h,p}(\mathbb{R}),\; i=1,2$ ({\bf commutative});
				%  \item $X_{h,p}(\mathbb{R})$ contains a zero element $\tilde 0$ such that 
				%  \begin{equation}
				% 	\tilde 0\otimes \tilde{x}=\tilde{x}\otimes \tilde 0 =\tilde 0  \quad\hbox{for all }\; \tilde{x}\in X_{h,p}(\mathbb{R}).
				%  \end{equation}
				%  Moreover, $\tilde 0=\langle 0; 1,1,0,0\rangle$.
				 \item $X_{h,p}(\mathbb{R})$ contains an identity element $\tilde{1}_X=\langle 1; e,e,1,1\rangle$ such that 
				 \begin{equation}\label{tilde1}
					\tilde{1}_X\otimes \tilde{x}=\tilde{x}\otimes \tilde{1}_X  \quad\hbox{for all }\; \tilde{x}\in X_{h,p}(\mathbb{R}).
				 \end{equation}
			\end{enumerate}
		\end{theorem}
% \begin{remark}
% 	1. \red{explain the zero element}
% \end{remark}

{\bf Proof:} Set $\tilde{x}_{i}=\langle x_i; d_i^-,d_i^+,\mu_i^-,\mu_i^+\rangle, i=1,2,3$. We give the proof one by one.
\begin{enumerate}[(i)]
    \item This is Theorem \ref{abeliangroup}.
    \item  Recall Formula $(4)$ of Definition \ref{def}, it holds
	\begin{equation*}
		\begin{array}{lll}
			(\tilde{x}_{1}\otimes\tilde{x}_{2})\otimes\tilde{x}_{3} 
			& = & \langle x_1x_2; (d_1^-)^{\ln d_2^-},(d_1^+)^{\ln d_2^+},\mu_1^-\mu_2^-,\mu_1^+\mu_2^+\rangle \otimes \langle x_3; d_3^-,d_3^+,\mu_3^-,\mu_3^+\rangle\\
			& = & \left\langle x_1x_2x_3; (d_1^-)^{\ln d_2^-\ln d_3^-},(d_1^+)^{\ln d_2^+\ln d_3^+},\mu_1^-\mu_2^-\mu_3^-,\mu_1^+\mu_2^+\mu_3^+\right\rangle \\
            & = & \left\langle x_1x_2x_3; (d_1^-)^{\ln [(d_2^-)^{\ln d_3^-}]}, (d_1^+)^{\ln [(d_2^+)^{\ln d_3^+}]},\mu_1^-\mu_2^-\mu_3^-,\mu_1^+\mu_2^+\mu_3^+\right\rangle \\
			& = & \tilde{x}_{1}\otimes(\tilde{x}_{2}\otimes\tilde{x}_{3}).
		\end{array}
    \end{equation*}
    \item Recall Formula $(1)$ and $(4)$ of Definition \ref{def}, we first compute the left distributive law: 
    \begin{equation*}
		\begin{array}{lll}
			&&\tilde{x}_1\otimes(\tilde{x}_2\oplus\tilde{x}_3)\\
            & = &
            \tilde{x}_1\otimes \langle x_2+x_3; d_2^-+d_3^-,d_2^++d_3^+,\mu_2^-+\mu_3^-,\mu_2^++\mu_3^+\rangle \\
            & = &
            \langle x_1(x_2+x_3); (d_1^-)^{\ln(d_2^-d_3^-)},(d_1^+)^{\ln(d_2^+d_3^+)},\mu_1^-(\mu_2^-+\mu_3^-),\mu_1^+(\mu_2^++\mu_3^+)\rangle \\
            & = &
            \langle x_1x_2+x_1x_3; (d_1^-)^{\ln d_2^-}(d_1^-)^{\ln d_3^-},(d_1^+)^{\ln d_2^+}(d_1^+)^{\ln d_3^+},
            \mu_1^-\mu_2^-+\mu_1^-\mu_3^-,\mu_1^+\mu_2^++\mu_1^+\mu_3^+\rangle \\
            & = &
            \tilde{x}_1\otimes\tilde{x}_2\oplus\tilde{x}_1\otimes\tilde{x}_3.
		\end{array}
    \end{equation*}
    The right distributive laws can be verified similarly.
    \item Formula $(4)$ of Definition \ref{def} gives
    \begin{equation*}
        \begin{array}{lll}
			\tilde{x}_{1}\otimes\tilde{x}_{2}
			& = & \langle x_1x_2; (d_1^-)^{\ln d_2^-},(d_1^+)^{\ln d_2^+},\mu_1^-\mu_2^-,\mu_1^+\mu_2^+\rangle \\
            & = & \langle x_2x_1; (d_2^-)^{\ln d_1^-},(d_2^+)^{\ln d_1^+},\mu_2^-\mu_1^-,\mu_2^+\mu_1^+\rangle \\
            & = & \tilde{x}_{2}\otimes\tilde{x}_{1}.
		\end{array}
    \end{equation*}
    Here we use the fact that $a^{\ln b}=b^{\ln a}$ for all $a,b>0$.
    \item For any $\tilde{x}=\langle x; d^-,d^+,\mu^-,\mu^+\rangle$, it holds that 
    $$
    \langle 1;e,e,1,1\rangle \otimes \langle x; d^-,d^+,\mu^-,\mu^+\rangle =\langle 1x; e^{\ln{d^-}},e^{\ln{d^+}},1\mu^-,1\mu^+\rangle = \langle x; d^-,d^+,\mu^-,\mu^+\rangle
    $$
    and the second identity holds since $a^{\ln e}=a$ for any $a\in \mathbb{R}^+$.
\end{enumerate}

Now we extend the binary operations on finite arithmetic operations  with the following notations.

The {\bf standard $n$ sum} $\displaystyle\bigoplus_{i=1}^n\tilde{x}_i $ of $\tilde{x}_1, \cdots, \tilde{x}_n\in X_{h,p}(\mathbb{R})$ is defined inductively by 
\begin{equation*}
		\bigoplus_{i=1}^1\tilde{x}_1= \tilde{x}_1; \,\; \hbox{and for }\;  n>1, \; \left(\bigoplus_{i=1}^{n-1}\tilde{x}_i\right) \oplus\tilde{x}_{n}
\end{equation*}
for any positive integer $n$. 

Consequently, the following properties of ring hold:

\begin{theorem}\label{ringproperty}
	$(X_{h,p}(\mathbb{R}),\oplus,\otimes)$ has the following properties: (Properties of ring)
	\begin{enumerate}[(i)]
		\item  $\tilde 0\otimes \tilde x=\tilde x\otimes \tilde 0=\tilde 0$ for all $\tilde x\in X_{h,p}(\mathbb{R})$;
		\item $(-\tilde{x}_1)\otimes\tilde{x}_2=\tilde{x}_1\otimes(-\tilde{x}_2)=-(\tilde{x}_1\otimes\tilde{x}_2)$ for all $\tilde x_i\in X_{h,p}(\mathbb{R}), i=1,2$;
		\item $(-\tilde{x}_1)\otimes(-\tilde{x}_2)=\tilde{x}_1\otimes\tilde{x}_2$ for all $\tilde x_i\in X_{h,p}(\mathbb{R}), i=1,2$;
		\item $(\lambda\tilde{x}_1)\otimes \tilde{x}_2=\tilde{x}_1\otimes(\lambda\tilde{x}_2)
        = \lambda(\tilde{x}_1\otimes \tilde{x}_2)$ for all $\lambda\in \mathbb{R}$ and $\tilde x_i\in X_{h,p}(\mathbb{R}), i=1,2$;
		\item $
		\displaystyle 
		\left(\bigoplus_{i=1}^n\tilde{x}_i\right)  
        \otimes
		\left(\bigoplus_{j=1}^m\tilde{y}_j\right)
		=
		\bigoplus_{i=1}^n \bigoplus_{j=1}^m \left(\tilde{x}_i\otimes \tilde{y}_j\right)
		$
		for all $\tilde{x}_i, \tilde{y}_j\in X_{h,p}(\mathbb{R})$, $i=1,\cdots, n$, $j=1,\cdots, m$.
	\end{enumerate}
\end{theorem}
{\bf Proof:} We give the proof one by one.
\begin{enumerate}[(i)]
    \item For any $\tilde{x}=\langle x; d^-,d^+,\mu^-,\mu^+\rangle$, it holds that 
    $$
    \langle 0; 1,1,0,0\rangle \otimes \langle x; d^-,d^+,\mu^-,\mu^+\rangle =\langle 0x; 1^{\ln{d^-}},1^{\ln{d^+}},0\mu^-,0\mu^+\rangle = \langle 0; 1,1,0,0\rangle 
    $$
    and the second identity holds since $a^{\ln 1}=1$ for any $a\in \mathbb{R}^+$.
    \item By the definition of $\otimes$, we compute
    \begin{equation*}
        \begin{array}{lll}
            (-\tilde{x}_1)\otimes\tilde{x}_2 
            &  =  &
            \langle - x_1; (d_1^-)^{-1},(d_1^+)^{-1}, -\mu_1^-, -\mu_1^+\rangle 
            \otimes
            \langle x_2;   d_2^-, d_2^+ ,\ \mu_2^-, \mu_2^+\rangle \\
            & = & \langle -x_1x_2; ((d_1^-)^{-1})^{{\ln d_2^-}},((d_1^+)^{-1})^{{\ln d_2^+}},-\mu_1^-\mu_2^-,-\mu_1^+\mu_2^+\rangle \\
            & = & \langle x_1(-x_2); (d_1^-)^{-\ln d_2^-},(d_1^+)^{-\ln d_2^+},\mu_1^-(-\mu_2^-),\mu_1^+(-\mu_2^+)\rangle \\
            & = & \tilde{x}_{1}\otimes(-\tilde{x}_{2})
        \end{array}
    \end{equation*}
    and  
    \begin{equation*}
        \begin{array}{lll}
            (-\tilde{x}_1)\otimes\tilde{x}_2 
            & = & \langle -(x_1x_2); ((d_1^-)^{\ln d_2^-})^{-1},((d_1^+)^{\ln d_2^+})^{-1},-(\mu_1^-\mu_2^-),-(\mu_1^+\mu_2^+)\rangle \\
            & = & -(\tilde{x}_1\otimes\tilde{x}_2).
        \end{array}
    \end{equation*}
    \item We compute
    \begin{equation*}
        \begin{array}{lll}
            (-\tilde{x}_1)\otimes(-\tilde{x}_2 )
            &  =  &
            \langle - x_1; (d_1^-)^{-1},(d_1^+)^{-1}, -\mu_1^-, -\mu_1^+\rangle 
            \otimes
            \langle- x_2; (d_2^-)^{-1},(d_2^+)^{-1}, -\mu_2^-, -\mu_2^+\rangle \\
            & = & \langle x_1x_2; ((d_1^-)^{-1})^{{\ln (d_2^-)^{-1}}},((d_1^+)^{-1})^{{\ln (d_2^+)^{-1}}},\mu_1^-\mu_2^-,\mu_1^+\mu_2^+\rangle \\
            & = & \langle x_1 x_2; (d_1^-)^{\ln d_2^-},(d_1^+)^{\ln d_2^+},\mu_1^-\mu_2^-,\mu_1^+\mu_2^+\rangle \\
            & = & \tilde{x}_{1}\otimes\tilde{x}_{2}.
        \end{array}
    \end{equation*}
    \item By the definition of scalar multiplication, we compute
    \begin{equation*}
        \begin{array}{lll}
            (\lambda\tilde{x}_1)\otimes\tilde{x}_2 
            &  =  &
            \langle - x_1; (d_1^-)^{\lambda},(d_1^+)^{\lambda}, \lambda\mu_1^-, \lambda\mu_1^+\rangle 
            \otimes
            \langle x_2;   d_2^-, d_2^+ ,\ \mu_2^-, \mu_2^+\rangle \\
            & = & \langle \lambda x_1x_2; ((d_1^-)^{\lambda})^{{\ln d_2^-}},((d_1^+)^{\lambda})^{{\ln d_2^+}},\lambda\mu_1^-\mu_2^-,\lambda\mu_1^+\mu_2^+\rangle \\
            & = & \langle x_1(\lambda x_2); (d_1^-)^{\lambda \ln d_2^-},(d_1^+)^{\lambda \ln d_2^+},\mu_1^-(\lambda \mu_2^-),\mu_1^+(\lambda \mu_2^+)\rangle \\
            & = & \tilde{x}_{1}\otimes(\lambda \tilde{x}_{2})
        \end{array}
    \end{equation*}
    and 
    \begin{equation*}
        \begin{array}{lll}
            (\lambda\tilde{x}_1)\otimes\tilde{x}_2 
            & = & \langle \lambda (x_1 x_2); \left((d_1^-)^{\ln d_2^-}\right)^\lambda,\left((d_1^+)^{\ln d_2^+}\right)^\lambda,\lambda (\mu_1^-\mu_2^-),\lambda (\mu_1^+  \mu_2^+)\rangle \\
            & = & \lambda(\tilde{x}_1\otimes \tilde{x}_2).
        \end{array}
    \end{equation*}
    \item We prove it by induction. The assertion is trivial for $m=n=1$. Assume that 
    \begin{equation*}
        (\tilde{x}_{1}\oplus \cdots\oplus \tilde{x}_{n})\otimes (\tilde{y}_{1}\oplus \cdots\oplus \tilde{y}_{m})
        =
        \bigoplus_{i=1}^n \bigoplus_{j=1}^m \left(\tilde{x}_i\otimes \tilde{y}_j\right)
    \end{equation*}
    holds for all  $n\leq N$ and $m\leq M$. Then by the distributive law (Formula $(iii)$ of Theorem \ref{ring}), it holds
    \begin{equation*}
        \begin{array}{lll}
            &  & 
            (\tilde{x}_{1}\oplus \cdots\oplus \tilde{x}_{N}\oplus \tilde{x}_{N+1})\otimes (\tilde{y}_{1}\oplus \cdots\oplus \tilde{y}_{M})\\
            & = &
             (\tilde{x}_{1}\oplus \cdots\oplus \tilde{x}_{N})\otimes (\tilde{y}_{1}\oplus \cdots\oplus \tilde{y}_{M})
            \oplus \tilde{x}_{M+1}\otimes (\tilde{y}_{1}\oplus \cdots\oplus \tilde{y}_{N})\\
            & = &
            \displaystyle\bigoplus_{i=1}^N \bigoplus_{j=1}^M \left(\tilde{x}_i\otimes \tilde{y}_j\right)
            \oplus (\tilde{x}_{N+1}\otimes \tilde{y}_{1}) \oplus\cdots  \oplus   (\tilde{x}_{N+1}\otimes \tilde{y}_{M})\\
            & = &
            \displaystyle\bigoplus_{i=1}^{N+1} \bigoplus_{j=1}^{M} \left(\tilde{x}_i\otimes \tilde{y}_j\right).
        \end{array}
    \end{equation*}
    Similarly, 
    \begin{equation*}
            (\tilde{x}_{1}\oplus \cdots\oplus \tilde{x}_{N})\otimes (\tilde{y}_{1}\oplus \cdots\oplus \tilde{y}_{M}\oplus \tilde{y}_{M+1})
        =
            \bigoplus_{i=1}^{N} \bigoplus_{j=1}^{M+1} \left(\tilde{x}_i\otimes \tilde{y}_j\right).
    \end{equation*}
    Hence the assertion holds true for any $m,n\in\mathbb{N}$.
    \endpf 
\end{enumerate}

We now consider those elements in $X_{h,p}(\mathbb{R})$ which have an inverse under multiplication. They form a group, and this  ``group of units'' is very important in algebraic structure theory. We start with the following definitions:
	\begin{definition}\label{invert}
			An element $\tilde{x}$ in $X_{h,p}(\mathbb{R})$ with is said to be {\bf left}{ (resp. right)} {\bf invertible} if there exists $\tilde{y}\in X_{h,p}(\mathbb{R}) $ (resp. $\tilde{z}\in X_{h,p}(\mathbb{R})$) such that $\tilde{y}\otimes \tilde{x}= \tilde{1}_X$ (reps. $\tilde{x} \otimes \tilde{z}=  \tilde{1}_X$). The element $\tilde{y}$ (resp. $\tilde{z}$) is called a {\bf left} (resp. right) {\bf inverse}  of $\tilde{x}$. An element $\tilde{x}\in X_{h,p}(\mathbb{R})$ that is both left and right invertible is denoted by $\tilde{x}^{-1}$ and said to be {\bf invertible} or to be a {\bf unit}. 
	\end{definition}
	
	\begin{theorem}\label{SpaceUX}
			Let 
			\begin{equation}\label{ux}
					U(X)=\{\tilde{x}=\langle x; d^-,d^+,\mu^-,\mu^+\rangle\in X_{h,p}(\mathbb{R})  \; |\;  x\neq 0,  d^-\neq 1,d^+\neq 1,\mu^-\neq 0,\mu^+\neq 0 \}.
			\end{equation}
			Then it is the set of all units in $X_{h,p}(\mathbb{R})$. Moreover,  for any $\tilde{x}\in U(X)$, 
			\begin{equation}\label{cancellation}
					\tilde{x}\otimes\tilde{y} = \tilde{0} \; \hbox{or } \;
					\tilde{y}\otimes\tilde{x} = \tilde{0}
					\;\Longrightarrow\;
					\tilde{y} = \tilde{0}.
			\end{equation}
	\end{theorem}

	\begin{remark}\label{nocancellation}
		Note that for $\tilde{x}_{1},\tilde{x}_{2}\in X_{h,p}(\mathbb{R})$, the following assertion is false:
		$$%\begin{equation}
	 \tilde{x}_1\otimes\tilde{x}_2 = \tilde{0}
	 \;\Longrightarrow\;
	 \tilde{x}_1 = \tilde{0}\quad \hbox{or } \quad \tilde{x}_2 = \tilde{0}.
			$$%\end{equation}
			In fact, any $\tilde{x}_{2}=\langle 0; d_2^-,d_2^+,\mu_2^-,\mu_2^+\rangle$, with the leading term  $x_2=0$,  acts as both a left and right zero divisor of the nonzero element $\tilde{x}_{1}=\langle x_1; 1,1,0,0\rangle$. According to the definition of the zero element $\tilde{0}$ in $X_{h,p}(\mathbb{R})$, the following two expressions are equivalent:
			$$
			\tilde{x}_1\otimes\tilde{x}_2 = \tilde{0}
		 \Longleftrightarrow 
			x_1x_2=0,\;\ln d_1^-\ln d_2^- = 0,\; \ln d_1^+\ln d_2^+ = 0,\;	\mu_1^-\mu_2^-     =0,\; 	\mu_1^+\mu_2^+ = 0.
			$$
			Consequently, all elements in $U(X)$ are invertible.
		\end{remark}
		\begin{remark}
			Note that the right and left cancellation laws hold in $U(X)$ by Assertion \eqref{cancellation}, that is, for all $\tilde{x}_i\in U(X), i=1,2,3$, 
			$$
			\tilde{x}_1\otimes \tilde{x}_2=\tilde{x}_1\otimes \tilde{x}_3
			\quad 	\hbox{or}	\quad
			\tilde{x}_2\otimes \tilde{x}_1=\tilde{x}_3\otimes \tilde{x}_1
			\quad\Longrightarrow\quad
			\tilde{x}_2=\tilde{x}_3.
			$$
		\end{remark}
		
	{\bf Proof of Theorem \ref{SpaceUX}:} For any $\tilde{x}\in U(X)$, taking into account that $\tilde{1}_X=\langle 1; e,e,1,1\rangle$, the inverse of $\tilde{x}$ is uniquely given by 
	\begin{equation*}
			\tilde{x}^{-1}=\langle \frac{1}{x}; \exp(\frac{1}{\ln d^-}),\exp(\frac{1}{\ln d^+}),\frac{1}{\mu^-},\frac{1}{\mu^+}\rangle
	\end{equation*}
	Clearly, \eqref{cancellation} is trivial. 
	\endpf

	As a direct consequence of Theorem \ref{SpaceUX}, it is possible to update the result for the first order fuzzy linear equation \eqref{1order} in Proposition \ref{linearequation1}. We have
	\begin{proposition}\label{linearequation2}
		Let $\tilde{a}=\langle a; d_a^-,d_a^+,\mu_a^-,\mu_a^+\rangle$ be in $U(X)$, $\tilde{b}=\langle b; d_b^-,d_b^+,\mu_b^-,\mu_b^+\rangle$ and $\tilde{c}=\langle c; d_c^-,d_c^+,\mu_c^-,\mu_c^+\rangle$ be in $X_{h,p}(\mathbb{R})$.  Set $\tilde{x}=\langle x; d^-,d^+,\mu^-,\mu^+\rangle$ be an undetermined fuzzy number such that 
	\begin{equation}\label{1order2}
		\tilde{a}\tilde{x}\oplus\tilde{b}=\tilde{c}.
	\end{equation}
	Then there exists a unique fuzzy number $\tilde{x}$ belonging to $X_{h,p}(\mathbb{R})$ satisfying \eqref{1order2}. More precisely, 
	\begin{equation}\label{Soln2}
		\tilde{x}%=\tilde{a}^{-1}\otimes(\tilde{c}\ominus\tilde{b})
			=
		\left\langle \frac{c-b}{a};  \exp\left(\frac{\ln d_c^--\ln d_b^-}{\ln d_a^-}\right) ,\exp\left(\frac{\ln d_c^+-\ln d_b^+}{\ln d_a^+}\right),
		\frac{\mu_c^--\mu_b^-}{\mu_a^-},\frac{\mu_c^--\mu_b^-}{\mu_a^+}\right\rangle.
	\end{equation}
	\end{proposition}

	The next binomial theorem is frequently useful in computations.  We define $\tilde{x}^0$  be the identity element $\tilde{1}_X$. The element $\tilde{x}^{n}$ is defined to be the standard $n$ product $\displaystyle\tilde{x}\otimes \tilde{x}\otimes\cdots\otimes \tilde{x}$.
	Recall that if $k$ and $n$ are integers with $0\leq k\leq n$, then the {\bf binomial coefficient } $\displaystyle\binom{n}{k}$ is the number $n!/(n-k)!k!$, where $0!=1$ and $n!=n(n-1)(n-2)\cdots 2\cdot 1$ for $n\geq 1$. It is obvious that $\displaystyle\binom{n}{k}$ is actually an integer.

	\begin{theorem}\label{Binomial}
		{\bf{(Binomial Theorem)}}.
		Let $n$ be a positive integer and $\tilde{x},\tilde{y},\tilde{x}_1,\tilde{x}_2,\cdots, \tilde{x}_n\in X_{h,p}(\mathbb{R})$.
		\begin{enumerate}[(i)]
			\item  %If $\tilde{x}\otimes\tilde{y}=\tilde{y}\otimes\tilde{x}$
			$
			\displaystyle
			(\tilde{x}\oplus\tilde{y})^n=\bigoplus_{k=0}^{n}\binom{n}{k}  \tilde{x}^k\otimes \tilde{y}^{n-k}
			$;
			\item Moreover,
			\begin{equation*}
				(\tilde{x}_1\oplus\cdots\oplus\tilde{x}_s)^n=\bigoplus\frac{n!}{(i_1!)\cdots(i_s!)}  \tilde{x}_1^{i_1}\otimes \tilde{x}_2^{i_2}\cdots \otimes\tilde{x}_s^{i_s},
			\end{equation*}
			where the sum is over all s-tuples $(i_1, i_2, \cdots, i_s)$ such that 
			$i_1+i_2+\cdots+i_s=n$.
		\end{enumerate}
	\end{theorem}
	{\bf Proof:} 	Assertion $(i)$ can be proved using induction on $n$ and the fact that $\tilde{x}$ and $\tilde{y}$ satisfy the distributive law and the commutativity. In fact, we compute
	\begin{equation*}
		\begin{array}{lll}
				\displaystyle (\tilde{x}\oplus\tilde{y})^{n+1}
				& = &
				\displaystyle (\tilde{x}\oplus\tilde{y})^n\otimes (\tilde{x}\oplus\tilde{y})\\
				& = &
				\displaystyle  \left(\bigoplus_{k=0}^{n}\binom{n}{k}  \tilde{x}^k\otimes \tilde{y}^{n-k}\otimes \tilde{x}\right)\oplus\left(\bigoplus_{k=0}^{n}\binom{n}{k}  \tilde{x}^k\otimes \tilde{y}^{n-k}\otimes \tilde{y}\right)\\
				& = &
				\displaystyle \tilde{x}^{n+1}\oplus\left(\bigoplus_{k=0}^{n-1}\binom{n}{k}  \tilde{x}^{k+1}\otimes \tilde{y}^{n-k}\right)\oplus\left(\bigoplus_{k=1}^{n}\binom{n}{k}  \tilde{x}^k\otimes \tilde{y}^{n-k+1}\oplus \tilde{y}^{n+1}\right)\oplus\tilde{y}^{n+1}\\
				& = &
				\displaystyle \tilde{x}^{n+1}\oplus\left(\bigoplus_{k=0}^{n-1}\left(\binom{n}{k}+\binom{n}{k+1}\right)  \tilde{x}^{k+1}\otimes \tilde{y}^{n-k}\right)\oplus\tilde{y}^{n+1}\\
				& = &
				\displaystyle \tilde{x}^{n+1}\oplus\left(\bigoplus_{k=0}^{n-1}\binom{n+1}{k+1} \tilde{x}^{k+1}\otimes \tilde{y}^{n-k}\right)\oplus\tilde{y}^{n+1}\\
				& = &
				\displaystyle \tilde{x}^{n+1}\oplus\left(\bigoplus_{k=1}^{n}\binom{n+1}{k} \tilde{x}^{k}\otimes \tilde{y}^{n+1-k}\right)\oplus\tilde{y}^{n+1}\\
				& = & \displaystyle \bigoplus_{k=0}^{n+1}\binom{n+1}{k} \tilde{x}^{k}\otimes \tilde{y}^{n+1-k}.
		\end{array}
\end{equation*}
Note that we use the formula that for $k<n$, it holds
	\begin{equation*}
			\begin{array}{lll}
					\displaystyle \binom{n}{k}+\binom{n}{k+1}
					& = &
					\displaystyle \frac{n!}{(n-k)!k!}+\frac{n!}{(n-k-1)!(k+1)!}\\
					& = &
					\displaystyle \frac{n!(k+1)}{(n-k)!(k+1)!}+\frac{n!(n-k)}{(n-k)!(k+1)!}\\
					& = &
					\displaystyle \frac{n!(n+1)}{(n-k)!(k+1)!}%\\
					%& = &
					= \displaystyle \binom{n+1}{k+1}.
			\end{array}
	\end{equation*}
	Assertion $(ii)$ can be proved using induction on $s$. The case $s=2$ is just part $(i)$ since
	\begin{equation*}
			(\tilde{x}_1\oplus\tilde{x}_2)^n=\bigoplus_{k=0}^{n}\binom{n}{k}  \tilde{x}_1^k\otimes \tilde{x}_2^{n-k}
			=
			\bigoplus_{k+j=n} \frac{n!}{k! j!}\tilde{x}_1^k\otimes \tilde{x}_2^{j}.
	\end{equation*}
	Assume that $(ii)$ is true for $s$. Note that
	\begin{equation*}
			\begin{array}{lll}
					(\tilde{x}_1\oplus\cdots\oplus\tilde{x}_s\oplus\tilde{x}_{s+1})^n
					&  =  &
					((\tilde{x}_1\oplus\cdots\oplus\tilde{x}_s)\oplus\tilde{x}_{s+1})^n\\
					&  =  &
					\displaystyle  \bigoplus_{k=0}^{n}\binom{n}{k} (\tilde{x}_1\oplus\cdots\oplus\tilde{x}_s)^k\otimes \tilde{x}_{s+1}^{n-k}\\
					&  =  &
					\displaystyle  \bigoplus_{k+j=n}\frac{n!}{k! j!}(\tilde{x}_1\oplus\cdots\oplus\tilde{x}_s)^k\otimes \tilde{x}_{s+1}^{j}\\
					&  =  &
					\displaystyle \bigoplus_{k+j=n}\frac{n!}{k! j!}\left(\bigoplus_{i_1+i_2+\cdots+i_s=k}\frac{k!}{(i_1!)\cdots(i_s!)}  \tilde{x}_1^{i_1}\otimes \tilde{x}_2^{i_2}\cdots \otimes\tilde{x}_s^{i_s}\right)\otimes \tilde{x}_{s+1}^{j}\\
					&  =  &
					\displaystyle \bigoplus_{i_1+i_2+\cdots+i_s+j=n} \frac{n!}{(i_1!)\cdots(i_s!)(j!)}  \tilde{x}_1^{i_1}\otimes \tilde{x}_2^{i_2}\cdots \otimes\tilde{x}_s^{i_s}\otimes \tilde{x}_{s+1}^{j}
			\end{array}
	\end{equation*}
	by part $(i)$. The proof is complete. \endpf
	
	As a simple consequence of Theorem \ref{Binomial} with $n=2$, we solve the following   fuzzy quadratic equation.
	\begin{proposition}\label{probc}
	Let $\tilde{b},\tilde{c}\in X_{h,p}(\mathbb{R})$ with $\tilde{b}=\langle b; d_b^-,d_b^+,\mu_b^-,\mu_b^+\rangle$ and $\tilde{c}=\langle c; d_c^-,d_c^+,\mu_c^-,\mu_c^+\rangle$. Set $\tilde{x}=\langle x; d^-,d^+,\mu^-,\mu^+\rangle$ be an undetermined fuzzy number such that 
	\begin{equation}\label{2order}
	\tilde{x}^2\ominus(\tilde{b}\oplus\tilde{c})\otimes\tilde{x}\oplus(\tilde{b}\otimes\tilde{c})=\tilde{0}.
	\end{equation}
	Then  $\tilde{x}=\tilde{b}$ and $\tilde{x}=\tilde{c}$ are solutions of the equation \eqref{2order}. 
	\end{proposition}
	\begin{remark}
			Note that the solvability of fuzzy quadratic equation is much more complicated than the standard quadratic equation due to the lack of cancellation law in $X_{h,p}(\mathbb{R})$ (see Remark \eqref{nocancellation}). Nevertheless, under some suitable conditions on the coefficients $ \tilde{a},\tilde{b},\tilde{c}$, one can give an quadratic formula for the fuzzy quadratic polynomial 
	\begin{equation}\label{quadratic}
			\tilde{a}\tilde{x}^2\oplus\tilde{b}\tilde{x}\oplus\tilde{c}=\tilde{0}.
	\end{equation}
	One possible assumption is that $\tilde{a}\in U(X)$ and Equation \eqref{quadratic} can be expressed as a product of 
	$$
	(\tilde{p}\tilde{x}\oplus\tilde{q})\otimes(\tilde{r}\tilde{x}\oplus\tilde{s})=\tilde{0}.
	$$
	Moreover, one can expect that the roots $\tilde{x}_{1}, \tilde{x}_{2}$ of the fuzzy quadratic polynomial $P(\tilde{x})= \tilde{a}\tilde{x}^2\oplus\tilde{b}\tilde{x}\oplus\tilde{c}$ satisfy 
	$$
	\tilde{x}_{1} \oplus \tilde{x}_{2}=\ominus \tilde{a}^{-1} \tilde{b}, \quad \tilde{x}_{1} \otimes\tilde{x}_{2}=  \tilde{a}^{-1} \tilde{c},
	$$
	which is the {\bf fuzzy Vieta's formulas } for $\tilde{a}\in U(X)$. Compare to the ``crisp" polynomial case, the main difference is that the fuzzy numbers are quinary arrays and a systematic analysis of the algebraic structure is needed.
	\end{remark}

\section{Vector space}\label{sec4}

We analyze the Gaussian-PDMF space with two operations: {\bf addition and scalar multiplication}. We claim that, 

\begin{theorem}\label{vectorspace}
	$X_{h,p}(\mathbb{R})$ is a vector space satisfying the following ten properties (axioms):
\begin{enumerate}[(1)]
	\item The sum of $\tilde{x}_1, \tilde{x}_2$, denoted by $\tilde{x}_1\oplus\tilde{x}_2$, is in $X_{h,p}(\mathbb{R})$.
	\item $\tilde{x}_1\oplus\tilde{x}_2=\tilde{x}_2\oplus\tilde{x}_1$.
	\item $(\tilde{x}_1\oplus\tilde{x}_2)\oplus\tilde{x}_3
			=
			\tilde{x}_1\oplus(\tilde{x}_2\oplus\tilde{x}_3)$.
	\item There is a {\bf zero} vector $\tilde{0}$ in $X_{h,p}(\mathbb{R})$ such that 
	$\tilde{x}\oplus\tilde{0}=\tilde{x}$ for any $\tilde{x}\in X_{h,p}(\mathbb{R})$.
	\item For each $\tilde{x}\in X_{h,p}(\mathbb{R})$, there is a vector $-\tilde{x}$ in $X_{h,p}(\mathbb{R})$ such that $\tilde{x}\oplus(-\tilde{x})=\tilde{0}$.
	\item The scalar multiplication of $\tilde{x}$ by $\lambda$, denoted by  $\lambda\tilde{x}$,  is in $X_{h,p}(\mathbb{R})$.
	\item $\lambda(\tilde{x}_1\oplus\tilde{x}_2)=\lambda\tilde{x}_1\oplus\lambda\tilde{x}_2$.
	\item $(\lambda_1+\lambda_2)\tilde{x}=\lambda_1\tilde{x}\oplus\lambda_2\tilde{x}$.
	\item $\lambda_1(\lambda_2\tilde{x})=(\lambda_1\lambda_2)\tilde{x}$.
	\item $1 \tilde{x}=\tilde{x}$.
\end{enumerate}
\end{theorem}
\begin{remark}
	Technically, $X_{h,p}(\mathbb{R})$ is a real vector space. All of the theory in this section should holds for a complex vector space if the scalars are complex numbers and Rule $(2)$ of Definition \ref{def} is adapted to complex numbers. The detail of the proof is beyond the scope of this paper and need to be done in the future.
\end{remark}

{\bf Proof of Theorem \ref{vectorspace}:} Assertion $(1)-(3), (6)$ are trivial. $(4)$ and $(5)$ is the assertion $(ii)$ and $(iii)$ of Theorem \ref{abeliangroup}.

For $(7)-(10)$, by Formula $(1)$ and $(2)$ of \eqref{def}, we compute
\begin{equation*}
    \begin{array}{lll}
        \lambda (\tilde{x}_{1}\oplus \tilde{x}_{2})
        & = &
        \lambda (\langle x_1; d_1^-,d_1^+,\mu_1^-,\mu_1^+\rangle\oplus \langle x_2; d_2^-,d_2^+,\mu_2^-,\mu_2^+\rangle)
        \\
        & = &
        \lambda (\langle x_1+x_2; d_1^-d_2^-,d_1^+d_2^+,\mu_1^-\mu_2^-,\mu_1^+\mu_2^+\rangle
        \\
        & = &
        \langle  \lambda ( x_1+x_2); (d_1^-d_2^-)^{\lambda},(d_1^+d_2^+)^{\lambda},\lambda(\mu_1^-+\mu_2^-),\lambda(\mu_1^++\mu_2^+)\rangle\\
        & = &
        \langle  \lambda   x_1+ \lambda x_2; (d_1^-)^{\lambda}(d_2^-)^{\lambda},(d_1^+)^{\lambda}(d_2^+)^{\lambda},\lambda\mu_1^-+\lambda\mu_2^-,\lambda\mu_1^++\lambda\mu_2^+)\rangle\\
        & = &
        \lambda\tilde{x}_1\oplus\lambda\tilde{x}_2,      
    \end{array}
\end{equation*}
and
\begin{equation*}
    \begin{array}{lll}
        (\lambda_1+\lambda_2 )\tilde{x}
        & = &
        \langle (\lambda_1+\lambda_2 ) x; (d^-)^{\lambda_1+\lambda_2 },(d^+)^{\lambda_1+\lambda_2 },(\lambda_1+\lambda_2 ) \mu^-,(\lambda_1+\lambda_2 ) \mu^+\rangle
        \\
        & = &
        \langle \lambda_1 x; (d^-)^{\lambda_1},(d^+)^{\lambda_1 }, \lambda_1\mu^-, \lambda_1 \mu^+ \rangle
        \oplus
        \langle  \lambda_2  x;  (d^-)^{\lambda_2 }, (d^+)^{ \lambda_2 },  \lambda_2  \mu^-,  \lambda_2  \mu^+\rangle
        \\
        & = &
        \lambda_1\tilde{x}\oplus\lambda_2\tilde{x},      
    \end{array}
\end{equation*}
and
\begin{equation*}
    \begin{array}{lll}
        \lambda_1(\lambda_2 \tilde{x})
        & = &
        \lambda_1\langle \lambda_2  x; (d^-)^{ \lambda_2 },(d^+)^{ \lambda_2 }, \lambda_2  \mu^-,\lambda_2   \mu^+\rangle
        \\
        & = &
        \langle \lambda_1\lambda_2  x; ((d^-)^{ \lambda_2 })^{\lambda_1},((d^+)^{ \lambda_2 })^{\lambda_1}, \lambda_1\lambda_2  \mu^-,\lambda_1\lambda_2   \mu^+\rangle
        \\
        & = &
        (\lambda_1 \lambda_2)\tilde{x},      
    \end{array}
\end{equation*}
and finally 
\begin{equation*}
    1 \tilde{x}= \langle 1 x; (d_1^-)^1,(d_1^+)^1,1 \mu_1^-,1 \mu_1^+\rangle = \tilde{x}.
\end{equation*}
\endpf 
In the sequel, we present a basis of the linear space $X_{h,p}(\mathbb{R})$. 

\begin{definition}
	An indexed set of vectors $\mathcal{X} = \{\tilde{x}_{1},\tilde{x}_{2},\cdots,\tilde{x}_{p}\}$ in $X_{h,p}(\mathbb{R})$ is said to be {\bf linearly independent } if the equation 
	\begin{equation*}
		\lambda_1 \tilde{x}_{1}\oplus \lambda_2 \tilde{x}_{2}\oplus\cdots\oplus\lambda_p \tilde{x}_{p}=\tilde{0}
	\end{equation*} 
	has only the trivial solution, $\lambda_1 =0, \cdots, \lambda_p=0$.
\end{definition}
Under the above definition, we have
\begin{theorem}\label{basis}
	A set of vectors $\mathcal{X} =\{\tilde{x}_{1},\tilde{x}_{2},\tilde{x}_{3},\tilde{x}_{4},\tilde{x}_{5}\}$ in $X_{h,p}(\mathbb{R})$ as the form 
    \begin{equation*}
    \tilde{x}_{1}=\langle 1;1,1,0,0\rangle,\; 
    \tilde{x}_{2}=\langle 0;e,1,0,0\rangle,\;
    \tilde{x}_{3}=\langle 0;1,e,0,0\rangle,\;
    \tilde{x}_{4}=\langle 0;1,1,1,0\rangle,\;
    \tilde{x}_{5}=\langle 0;1,1,0,1\rangle
    \end{equation*}
	%$$
	% \begin{array}{lll}
	% 	\tilde{x}_{1}=\langle 1;1,1,0,0\rangle\\
	% 	\tilde{x}_{2}=\langle 0;e,1,0,0\rangle\\
	% 	\tilde{x}_{3}=\langle 0;1,e,0,0\rangle\\
	% 	\tilde{x}_{4}=\langle 0;1,1,1,0\rangle\\
	% 	\tilde{x}_{5}=\langle 0;1,1,0,1\rangle\\
	% \end{array}
	% $$
	is a {\bf{basis}} for $X_{h,p}(\mathbb{R})$. More precisely, 
	\begin{enumerate}
		\item  $\mathcal{X} $ is a linear independent set, and
		\item $X_{h,p}(\mathbb{R})$ can be spanned by $\mathcal{X}$,i.e. for any $\tilde{x}\in X_{h,p}(\mathbb{R})$, there exists $\lambda_i\in \mathbb{R},i = 1,2,3,4,5$ such that 
		$$
		\tilde{x}=\sum_{i=1}^{5}\lambda_i \tilde{x}_{i}=
		\lambda_1 \tilde{x}_{1}\oplus
		\lambda_2 \tilde{x}_{2}\oplus
		\lambda_3 \tilde{x}_{3}\oplus
		\lambda_4 \tilde{x}_{4}\oplus
		\lambda_5 \tilde{x}_{5}.
		$$
	\end{enumerate}
	Moreover, the set of scalars $ \lambda_1,\cdots, \lambda_5$ is unique. 
\end{theorem}
We call that the {\bf coordinates of $\tilde{x}$ relative to the basis $\mathcal{X}$} (or the {\bf $\mathcal{X}$-coordinates of $\tilde{x}$}) are the weight  $\lambda_1, \cdots, \lambda_5$  such that 
$$
\tilde{x}=
\lambda_1 \tilde{x}_{1}\oplus
		\lambda_2 \tilde{x}_{2}\oplus
		\lambda_3 \tilde{x}_{3}\oplus
		\lambda_4 \tilde{x}_{4}\oplus
		\lambda_5 \tilde{x}_{5}.
$$

We call the vector in $\mathbb{R}^n$
	$$
	\[\tilde{x}\]_{\mathcal{X}} = 
	\[
			\lambda_1,
			\lambda_2,
			\lambda_3,
			\lambda_4,
			\lambda_5 
	\]
	$$
	the {\bf $\mathcal{X}$-coordinate vector of $\tilde{x}$}. The mapping $x \mapsto \[\tilde{x}\]_{\mathcal{X}} $ is the {\bf coordinate mapping determined by $\mathcal{X}$}. In fact, it connects the Gaussian-PDMF space $X_{h,p}(\mathbb{R})$ to the standard Euclidean space $\mathbb{R}^5$.  We have
	\begin{theorem}\label{xcoordinate}
		The $\mathcal{X}$-coordinate vector of $\tilde{x}=\langle x; d^-, d^+,\mu^-,\mu^+\rangle$ is given by 
		\begin{equation*}
			\[\tilde{x}\]_{\mathcal{X}}
			=
			\[x, \ln d^-, \ln d^+,\mu^-,\mu^+ \] 
		\end{equation*}
		for any $\tilde{x}\in X_{h,p}(\mathbb{R})$. Moreover, the coordinate mapping $x \mapsto \[\tilde{x}\]_{\mathcal{X}} $  is a one-to-one linear map (isomorphism) from $X_{h,p}(\mathbb{R})$ onto $\mathbb{R}^5$, i.e., for any $\lambda_1,\lambda_2\in\mathbb{R}$,
		\begin{equation*}
			\[\lambda_1\tilde{x}_{1}\oplus \lambda_2\tilde{x}_{2}\]_{\mathcal{X}}
			=
			\lambda_1\[\tilde{x}_{1}\]_{\mathcal{X}} + \lambda_2\[\tilde{x}_{2}\]_{\mathcal{X}}.
		\end{equation*}
	\end{theorem}
    {\bf Proof:} It is straightforward. 
    \endpf

		\section{Numerical simulation}\label{sec5}

		In this section, we illustrate three examples to show the advantage of the Gaussian-PDMF space with the algebraic structures. The first and the second example are linear first order fuzzy equation, with real numbers and fuzzy numbers coefficients, respectively. The third one is a quadratic monic fuzzy equation.

		{\bf Example $1$}. Let $\tilde{b},\tilde{c}\in X_{h,p}(\mathbb{R})$ with $\tilde{1}=\langle 1; 1,1,0,0\rangle$ and $\tilde{3}=\langle 3; 4,1,1,2\rangle$. Let $a=2$. Set $\tilde{x}=\langle x; d^-,d^+,\mu^-,\mu^+\rangle$ be an undetermined fuzzy number in $X_{h,p}(\mathbb{R})$ satisfying \eqref{1order}, i.e.,
		\begin{equation}\label{2x-1=3}
			2\tilde{x}\oplus\tilde{1}=\tilde{3}.
		\end{equation}
		
		{\bf Solution $A$}. Using Proposition \ref{linearequation1}, $\tilde{x}$ is given by 
		\begin{equation*}
			\tilde{x}=\frac{1}{2}(\tilde{3}\ominus\tilde{1})=
			\langle 1; 2,1,\frac{1}{2},1\rangle.
		\end{equation*}

\begin{figure}[htbp]
			\centering
			\begin{minipage}[t]{0.48\textwidth}
				\centering
				\includegraphics[width=6cm]{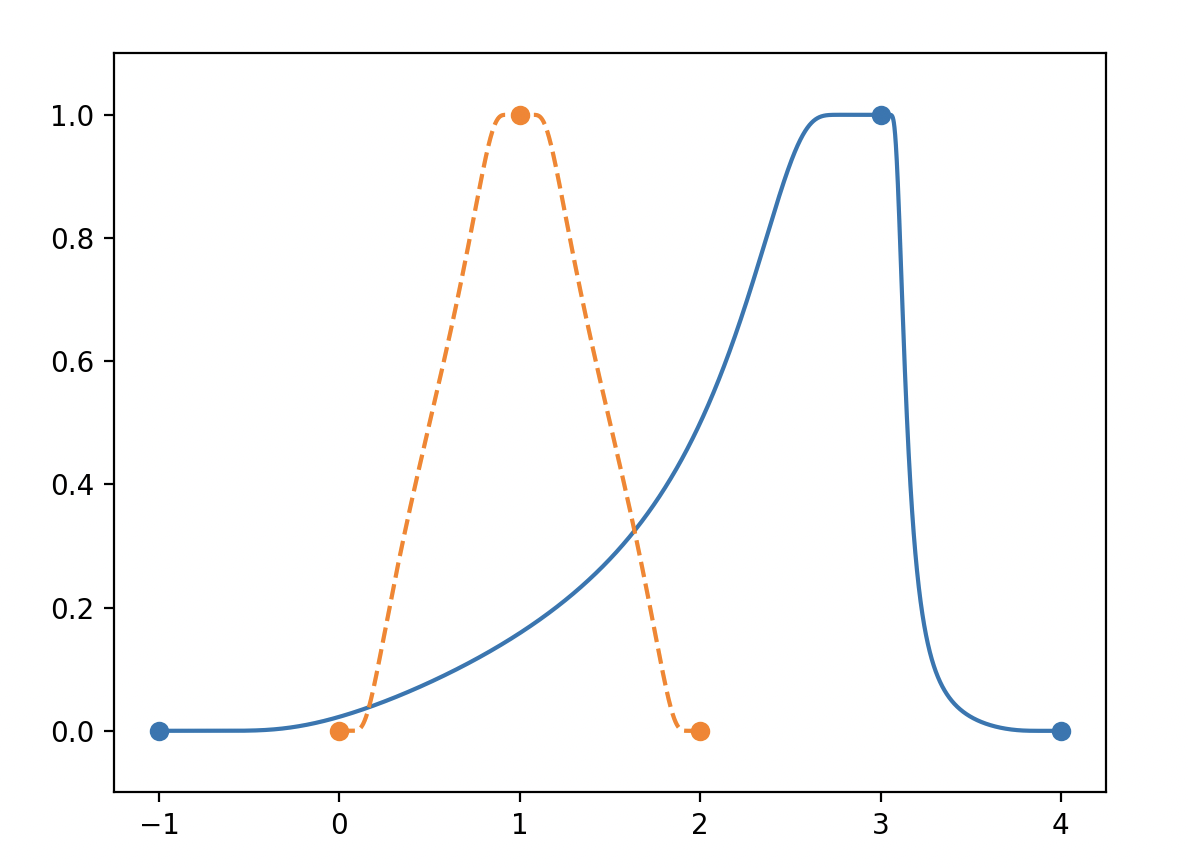}
				\caption{The graph of $\tilde{3}$ and $\tilde{1}$}\label{Fig.2} 
			\end{minipage}
			\begin{minipage}[t]{0.48\textwidth}
				\centering
				\includegraphics[width=6cm]{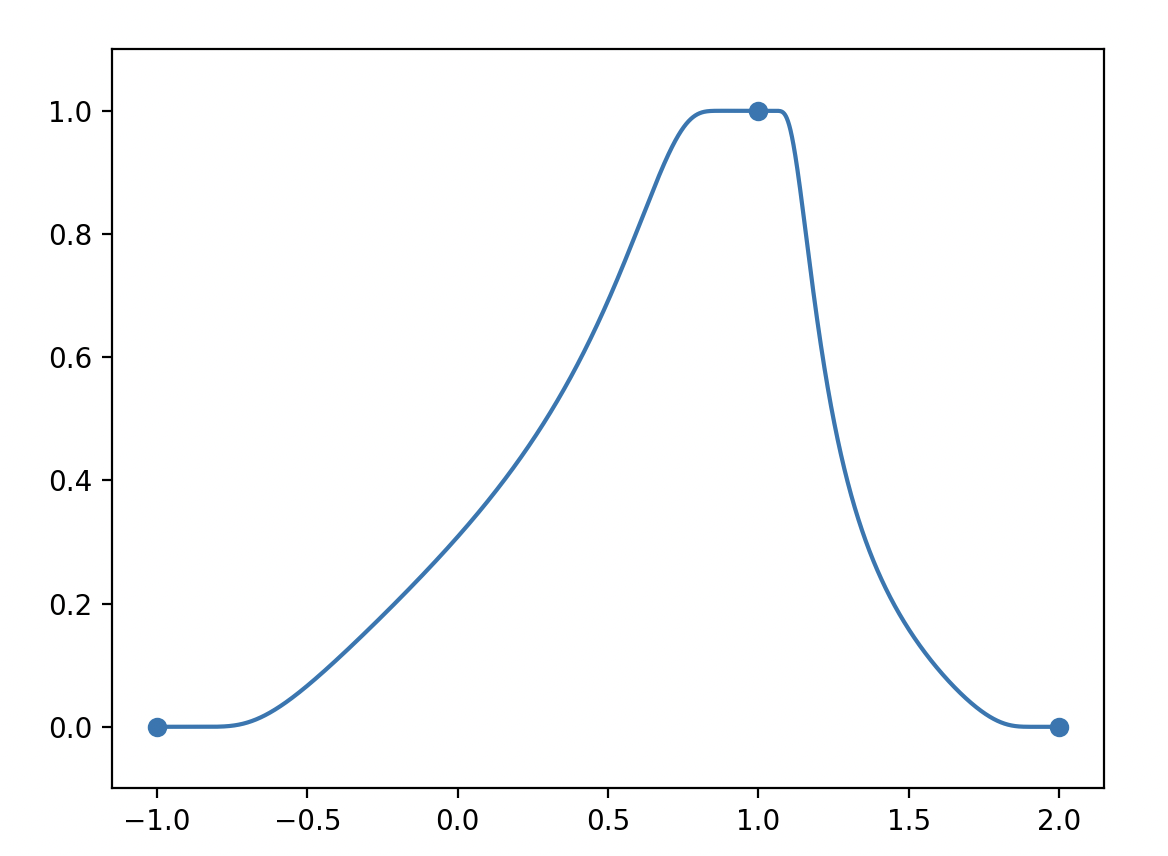}
				\caption{The graph of $\frac{1}{2}(\tilde{3}\ominus\tilde{1})$}\label{Fig.3} 
			\end{minipage}
		\end{figure}

		The corresponding fuzzy numbers are in Figure \ref{Fig.2} and Figure \ref{Fig.3}.

		{\bf Solution $B$}. Using Theorem \ref{xcoordinate}, we compute the $\mathcal{X}$-coordinate vector of $\tilde{x}$ as 
		$$
		\[\tilde{1}\]_{\mathcal{X}}=\[1,0,0,0,0\]_{\mathcal{X}}, \quad \[\tilde{3}\]_{\mathcal{X}}=\[3,\ln 4,0,1,2\]_{\mathcal{X}}
		$$
		\eqref{2x-1=3} holds if and only if 
		\begin{equation}\label{5dvector}
			2\[\tilde{x}\]_{\mathcal{X}}+ \[\tilde{1}\]_{\mathcal{X}}
					= 
			\[\tilde{3}\]_{\mathcal{X}}.
		\end{equation}
		\eqref{5dvector} is a standard arithmetic equation and can be solved as 
		$$
		\[\tilde{x}\]_{\mathcal{X}}=\[1, \ln ,0,\frac{1}{2},1\]_{\mathcal{X}},
		$$
		It means that
		\begin{equation*}
			\tilde{x}=	\langle 1; 2,1,\frac{1}{2},1\rangle.
		\end{equation*}
		
		We show an example of Proposition \ref{linearequation2}, which solves $\tilde{a}\tilde{x}\oplus\tilde{b}=\tilde{c}$.

\begin{figure}[htbp]
			\centering
			\begin{minipage}[t]{0.48\textwidth}
				\centering
				\includegraphics[width=6cm]{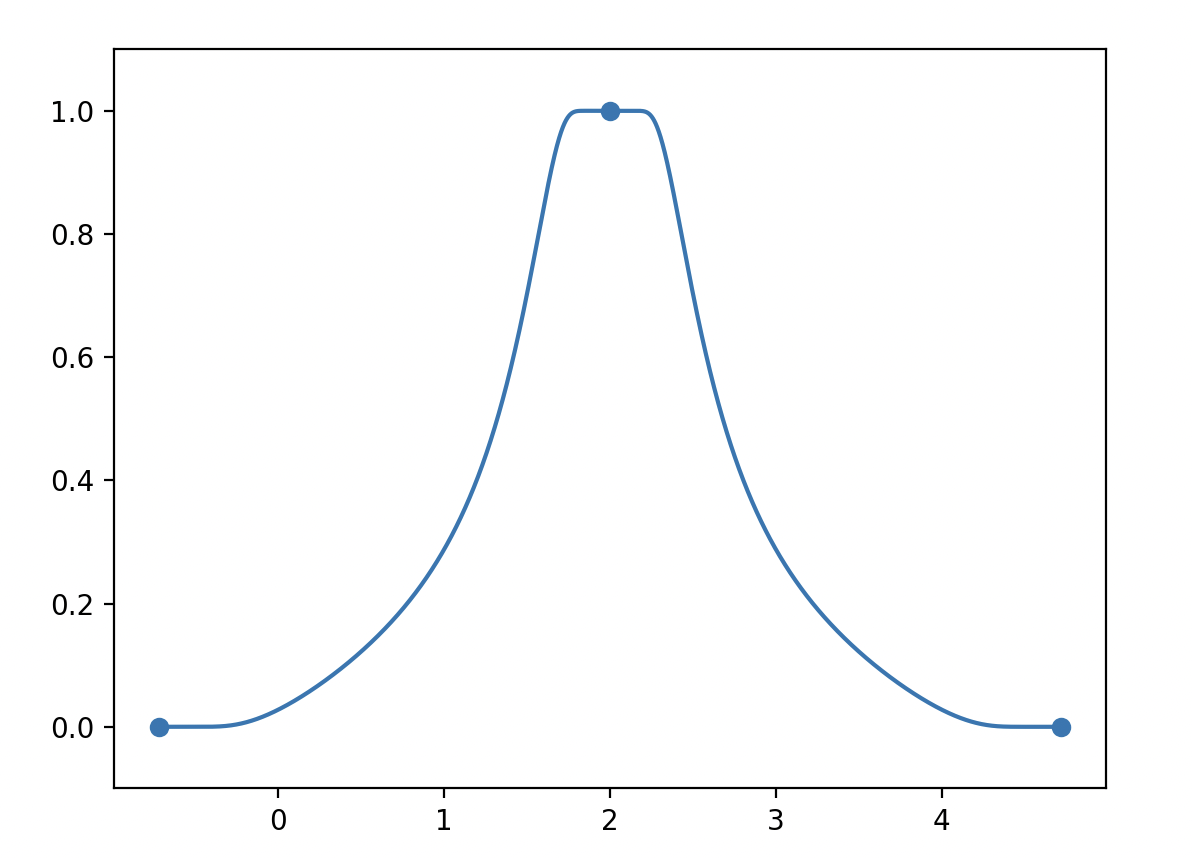}
				\caption{The graph of $\tilde{2}$}\label{Fig.4} 
			\end{minipage}
			\begin{minipage}[t]{0.48\textwidth}
				\centering
				\includegraphics[width=6cm]{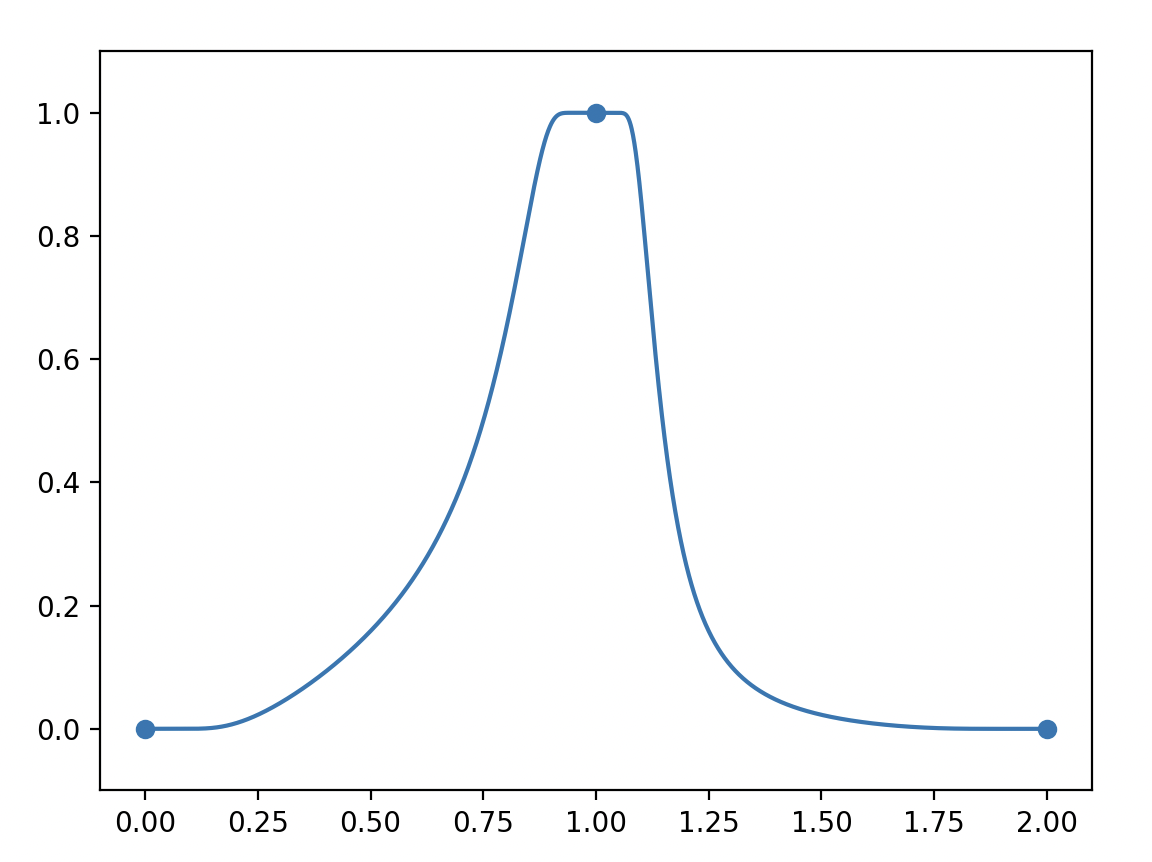}
				\caption{The graph of $\tilde{2}^{-1}\otimes (\tilde{3}\ominus\tilde{1})$}\label{Fig.5} 
			\end{minipage}
		\end{figure}

{\bf Example $2$}.  Let $\tilde{b},\tilde{c}\in X_{h,p}(\mathbb{R})$ with $\tilde{1}=\langle 1; 1,1,0,0\rangle$ and $\tilde{3}=\langle 3; 4,1,1,2\rangle$. Let $\tilde{a}=\langle 2; e,e,1,1\rangle$. Set $\tilde{x}=\langle x; d^-,d^+,\mu^-,\mu^+\rangle$ be an undetermined fuzzy number in $X_{h,p}(\mathbb{R})$ satisfying \eqref{1order2}, i.e.,
\begin{equation}\label{f2x-1=3}
	\tilde{2}\tilde{x}\oplus\tilde{1}=\tilde{3}.
\end{equation}
{\bf Solution}. Note that Theorem \ref{xcoordinate} is no longer applicable and we apply Proposition \ref{linearequation2}. We first compute 
$$
\tilde{2}^{-1}=\langle \frac12; 1,1,1,1\rangle, \quad \tilde{3}\ominus\tilde{1} =\langle 2; 4,1,1,2\rangle.
$$
Hence 
$$
\tilde{x}=\tilde{2}^{-1}\otimes (\tilde{3}\ominus\tilde{1})=\langle 1; 1,1,1,2\rangle.
$$
The corresponding fuzzy numbers are in Figure \ref{Fig.4} and Figure \ref{Fig.5}.

{\bf Example $3$}. Let $\tilde{2}=\langle 2; e,e,1,2\rangle$ and $\tilde{3}=\langle 3; 1,1,0,1\rangle$. Set $\tilde{x}=\langle x; d^-,d^+,\mu^-,\mu^+\rangle$ be an undetermined fuzzy number in $X_{h,p}(\mathbb{R})$ satisfying 
\begin{equation}\label{xsquare}
    \tilde{x}^2\oplus \tilde{2}\tilde{x}\ominus\tilde{3}=\tilde{0}.
\end{equation} 
{\bf Solution}. We solve the $5$-tuple $\tilde{x}$ one by one.
\begin{itemize}
    \item For $X_{h,p}(\mathbb{R})$, we have $x^2+2x-3=0$. Hence, $x_1=1$ or $x_2=-3$;
    \item For $d^-$, we have 
    \begin{equation*}
        (d^-)^{\ln d^-}\cdot e^{\ln d^-}\cdot 1^{-1}=1 \Leftrightarrow  
        (\ln d^-)^2+\ln d^-+0=0  \Leftrightarrow 
        \ln d^-(\ln d^-+1)=0.
    \end{equation*}
    Hence, $d^-_1=1$ or $d^-_2=e^{-1}$;
    \item For $d^+$, we have the same equation as $d^-$. Hence, $d^+_1=1$ or $d^+_2=e^{-1}$;
    \item For $\mu^-$, we have $(\mu^-)^2+\mu^-=0$. Hence, $\mu^-_1=-1$ or  $\mu^-_2=0$;
    \item For $\mu^+$, we have $(\mu^+)^2+2\mu^++1=0$. Hence, $\mu^+_{1,2}=-1$.
\end{itemize}
Interestingly, there are $2^4$ fuzzy numbers in $X_{h,p}(\mathbb{R})$ satisfying Equation \eqref{xsquare}:
\begin{equation*}
    \tilde{x}=\langle x_i; d_j^-,d_k^+,\mu_l^-,-1\rangle\;\hbox{for any } \; i,j,k,l=1,2.
\end{equation*} 

\section{Final remarks}\label{sec6}

In summary, our research introduces a new class of membership functions, merging subjective perception and objective information through the Gaussian Probability Density Membership Function (Gaussian-PDMF) space. This space, denoted as $X_{h,p}(\mathbb{R})$, uniquely identifies fuzzy numbers, symbolized by $\tilde{x}$, with a vector notation $\langle x; d^-, d^+, \mu^-,\mu^+\rangle$. We define five algebraic operations, revealing a robust structure within $X_{h,p}(\mathbb{R})$. Especially, addition $\oplus$ forms an abelian group, and multiplication $\otimes$ establishes a commutative ring with identity. This structured algebraic framework proves instrumental in solving fuzzy equations, applying the Binomial Theorem, and exploring fuzzy Vieta's formula.

Our work could extend theoretical insights into fuzzy numbers, providing a fresh framework for understanding and applying fuzzy mathematics. The practical utility of this algebraic structure is demonstrated through several numerical examples. We hope this research should contributes to the broader field of fuzzy mathematics and related research areas.

\bibliography{ref-fuzzy.bib}
\end{document}